\documentclass[11pt]{amsart}
\usepackage{latexcad}
\usepackage{natbib}

\def\ldiv{\setminus}
\def\rdiv{/}
\def\text{\textrm}

\theoremstyle{definition}
\newtheorem{definition}{Definition}[section]
\newtheorem{example}[definition]{Example}

\theoremstyle{plain}
\newtheorem{theorem}[definition]{Theorem}
\newtheorem{lemma}[definition]{Lemma}
\newtheorem{corollary}[definition]{Corollary}
\newtheorem{proposition}[definition]{Proposition}
\newtheorem{remark}[definition]{Remark}

\title[C-loops]{C-loops: An introduction}

\author{J.~D.~Phillips}

\address{Department of Mathematics \& Computer Science, Wabash College,
Crawfordsville, Indiana 47933, U.S.A.}

\email{phillipj@wabash.edu}

\author{Petr Vojt\v{e}chovsk\'y}

\address{Department of Mathematics, University of Denver, 2360 S Gaylord St,
Denver, CO, 80208, U.S.A.}

\email{petr@math.du.edu}

\begin{document}

\begin{abstract}
C-loops are loops satisfying $x(y(yz))=((xy)y)z$. They often behave
analogously to Moufang loops and they are closely related to Steiner triple
systems and combinatorics. We initiate the study of C-loops by proving: (i)
Steiner loops are C-loops, (ii) C-loops are alternative, inverse property
loops with squares in the nucleus, (iii) the nucleus of a C-loop is a normal
subgroup, (iv) C-loops modulo their nucleus are Steiner loops, (v) C-loops
are power associative, power alternative but not necessarily diassociative,
(vi) torsion commutative C-loops are products of torsion abelian groups and
torsion commutative $2$-C-loops; and several other results. We also give
examples of the smallest nonassociative C-loops, and explore the analogy
between commutative C-loops and commutative Moufang loops.
\end{abstract}



\maketitle

\section{Introduction}

\noindent \emph{C-loops} are loops satisfying the identity
\begin{equation}\label{Eq:C}
    x(y(yz)) = ((xy)y)z.
\end{equation}
As we shall see, they are in a sense dual to Moufang loops---the most
intensively studied variety of loops---and they are closely related to
Steiner triple systems. They are thus important both algebraically and
combinatorially, and they are amenable to analysis by techniques from both
fields. But in spite of this, little is known about them. It is the intention
of this paper to remedy this situation by laying a foundation for the
systematic study of C-loops.

We assume that the reader is familiar with the reasoning and notational
conventions of loop theory, however, we do not hesitate to include
loop-theoretical folklore and to point out some of the pitfalls of
nonassociativity---mostly because we fell into many of them ourselves.

C-loops were named by Ferenc Fenyves \cite{Fe2}, who investigated the
inclusions between varieties of loops of \emph{Bol-Moufang type}. These are
varieties of loops defined by a single identity that: (i) involves three
distinct variables on both sides, (ii) contains variables in the same order on
both sides, (iii) exactly one of the variables appears twice on both sides.

\begin{figure}[ht]\begin{center}\input{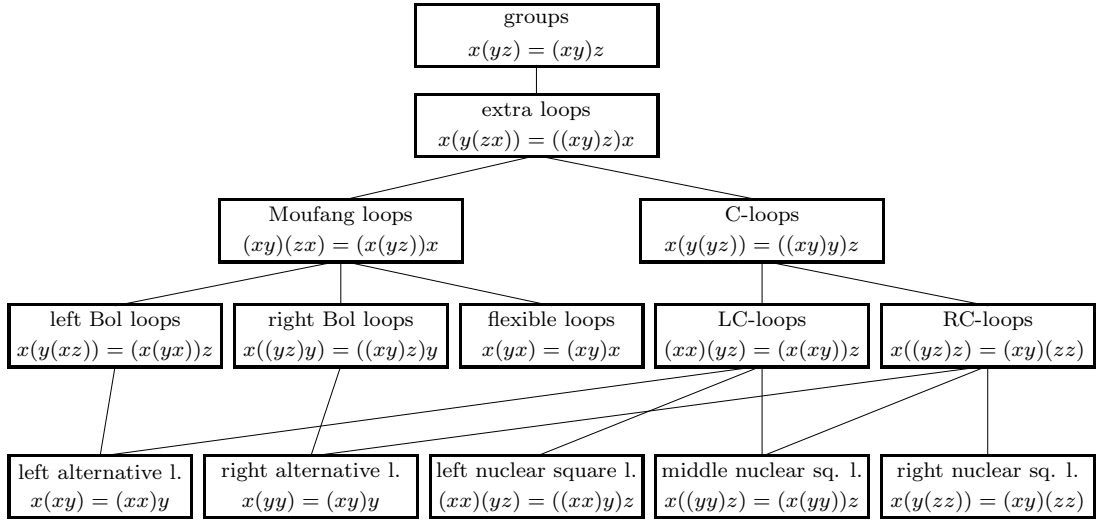}\end{center}
\caption{Varieties of loops of Bol-Moufang type.}\label{Fg:L}
\end{figure}

Fenyves's program was completed by the authors in
\cite{PhVo}. There are 60 identities of Bol-Moufang type,
and they happen to define 14 distinct varieties of loops.
Figure \ref{Fg:L} gives the Hasse diagram of these
varieties, with the largest varieties (with respect to
inclusion) at the bottom.

\setlength{\unitlength}{1mm}

A superficial glance at the diagram suggests that C-loops could behave
analogously to Moufang loops. This impression is further strengthened by the
fact that C-loops are exactly those loops that are both LC-loops and RC-loops
\cite[Theorem 4]{Fe2}, just as Moufang loops are exactly those loops that are
both left Bol and right Bol \cite{Bol}. There are additional analogies,
especially between commutative Moufang loops and commutative C-loops, as we
shall see.

\section{C-loops and Steiner loops}

\noindent In combinatorics, Moufang loops have connections
to projective geometry (Moufang planes, Moufang polygons,
etc., cf. \cite{Ma}), while C-loops have connections to
Steiner triple systems:

Consider the complete graph $K_n$ on $n$ vertices. A \emph{Steiner triple
system} is a decomposition of the edges of $K_n$ into disjoint triangles. It is
well known (cf.\ \cite{CR}) that such a decomposition exists if and only if
$n\equiv 1 \pmod 6$ or $n\equiv 3\pmod 6$; the case $n=1$ being degenerate.

There is a canonical way of constructing a quasigroup from a Steiner triple
system. Namely, if $S=K_n$ is a Steiner triple system, we define multiplication
on $\{1,\dots,n\}$ by $xx=x$, and (for $x\ne y$) by $xy=z$ if and only if
$\{x,y,z\}$ is a triangle of $S$. The resulting quasigroup clearly satisfies
\begin{equation}\label{Eq:SQ}
    xx=x,\quad (yx)x=y, \quad xy=yx.
\end{equation}
Conversely, any quasigroup satisfying $(\ref{Eq:SQ})$ gives rise to a Steiner
triple system in a canonical way (cf. \cite{CR}, \cite{Li}). Quasigroups
satisfying $(\ref{Eq:SQ})$ are therefore called \emph{Steiner quasigroups}.

Any Steiner quasigroup can be made into a loop by introducing a new element $e$
and by letting $xx=e$, $xe=ex=x$. Such loops satisfy
\begin{equation}\label{Eq:SL}
    xx=e,\quad (yx)x=y, \quad xy=yx,
\end{equation}
and are called \emph{Steiner loops}. It is now clear that the Steiner
quasigroup that gave rise to a Steiner loop $L$ can be reconstructed from $L$.
Steiner loops are therefore in one-to-one correspondence with Steiner triple
systems, too.

Intuitively, the reason why C-loops are related to Steiner loops is the
presence of the term $(xy)y$ in the defining equation $(\ref{Eq:C})$. More
formally:

\begin{lemma}\label{Lm:CS} Every Steiner loop is a C-loop.
\end{lemma}
\begin{proof}
Note that $(xy)y=x$ is a part of the definition $(\ref{Eq:SL})$, and that
$y(yz)=z$ follows from $(\ref{Eq:SL})$ immediately by commutativity. Thus
$x(y(yz))=xz=((xy)y)z$.
\end{proof}

Not every C-loop is a Steiner loop, as is witnessed by any nonabelian group.

Another connection between C-loops and Steiner loops
becomes apparent upon investigating the nucleus of C-loops.

Recall that for a loop $L$, the set $N_\lambda=\{x\in L;\;x(yz)=(xy)z$ for
every $y$, $z\in L\}$ is called the \emph{left nucleus}. Similarly, the
\emph{middle nucleus} $N_\mu$ consists of all elements $x\in L$ satisfying
$y(xz)=(yx)z$ for every $y$, $z\in L$; and the \emph{right nucleus} $N_\rho$
consists of all elements $x\in L$ satisfying $y(zx)=(yz)x$ for every $y$, $z\in
L$. The \emph{nucleus} $N=N_\lambda\cap N_\mu\cap N_\rho$ of $L$ is a subgroup
of $L$.

There are several equivalent ways in which normality can be defined for loops.
The following definition works best with elementary calculations. A subloop $K$
of a loop $L$ is said to be \emph{normal} in $L$ if $xK=Kx$, $x(yK)=(xy)K$, and
$x(Ky)=(xK)y$ for every $x$, $y\in L$. The \emph{factor loop} $L/K$ is then
defined in the usual way.

A loop $L$ with neutral element $e$ is a \emph{left inverse property loop} if
$x'(xy)=y$ for every $x$, $y\in L$, where $x'$ is the unique element satisfying
$x'x=e$. Dually, $L$ is a \emph{right inverse property loop} if $(yx)x''=y$ for
every $x$, $y\in L$, where $x''$ is the unique element satisfying $xx''=e$. A
loop that has both the left and right inverse property is an \emph{inverse
property loop}.

If $x\in L$ is such that $x'(xy)=(yx)x''=y$ for every $y$, we have
$x'=x'e=x'(xx'')=x''$. Therefore, inverse property loops possess
\emph{two-sided inverses} (i.e., $x'=x''=x^{-1}$), and it is easy to check that
they satisfy the \emph{antiautomorphic inverse property} (i.e.,
$(xy)^{-1}=y^{-1}x^{-1})$.

Pflugfelder shows \cite[p.\ 123]{Pf} that Steiner loops are exactly commutative
inverse property loops of exponent 2. In fact, Steiner loops are exactly
inverse property loops of exponent 2. This fact belongs to loop-theoretical
folklore and is sometimes used as a definition of Steiner loops (cf.
\cite{KKPARIF}). Since we did not manage to find a reference for the proof,
here it is:

\begin{lemma}\label{Lm:Folklore}
Steiner loops are exactly inverse property loops of exponent two.
\end{lemma}
\begin{proof} Let $L$ be a Steiner loop. Since $xx=e$, every element is its own
two-sided inverse. From $(yx)x=y$ we see that $L$ has the left inverse
property. By commutativity, it has the right inverse property, too.

Conversely, let $L$ be an inverse property loop of exponent $2$. Let $z=xy$.
Then $xz=x(xy)=x^{-1}(xy)=y$, and similarly, $x=yz$, $yx=z$. Thus $L$ is
commutative. As $(yx)x=y$ by the right inverse property, $L$ is a Steiner loop.
\end{proof}

Also notice that $xx=e$ is not necessary in the definition $(\ref{Eq:SL})$ of
Steiner loop, since $xx=(ex)x=e$. Hence Steiner loops are exactly loops
satisfying
\begin{equation}\label{Eq:TS}
    (yx)x=y,\quad xy=yx.
\end{equation}
Quasigroups satisfying $(\ref{Eq:TS})$ are called \emph{totally symmetric}, and
thus Steiner loops can also be found under the name \emph{totally symmetric
loops} in the literature.

Let us now mention some basic properties of LC-loops and
C-loops that we will use without reference throughout the
paper. The first three properties are due to Fenyves
\cite[Theorem 2]{Fe2}. The fourth property first appeared
in \cite{PhVo}.

\begin{proposition}\label{Pr:Sum} Let $L$ be an LC-loop. Then:
\begin{enumerate}
\item[(i)] $L$ is left alternative,

\item[(ii)] $L$ has the left inverse property,

\item[(iii)] $L$ is a left nuclear square loop,

\item[(iv)] $L$ is a middle nuclear square loop.
\end{enumerate}
\end{proposition}

We will often derive theorems from their one-sided versions.

\begin{corollary} Let $L$ be a C-loop. Then:
\begin{enumerate}
\item[(i)] $L$ is both left alternative and right alternative,

\item[(ii)] $L$ has the inverse property,

\item[(iii)] $L$ is a \emph{nuclear square loop}, i.e., $x^2$ belongs to the
nucleus of $L$ for every $x\in L$.
\end{enumerate}
\end{corollary}

\begin{corollary}\label{Cr:Nuc} The three nuclei of a C-loop coincide.
\end{corollary}
\begin{proof}
The three nuclei coincide for any inverse property loop, by \cite[Theorem
VII.2.1]{Br}.
\end{proof}

The nucleus $N$ of a loop $L$ is always a subgroup of $L$, but it is not
necessarily a normal subgroup of $L$. Even when $L$ is an inverse property
loop, its nucleus does not have to be normal in $L$. (See Example
\ref{Ex:IPNuc}).

Throughout the paper, if we claim without explanation that
a loop with given properties is as small as possible, or
that there are $m$ such nonisomorphic loops of given order,
we rely on the finite model builder Mace4 \cite{Mc}.

\begin{example}\label{Ex:IPNuc}
The smallest inverse property loop with nucleus that is not normal.
\begin{displaymath}
\begin{array}{cccccccccccc}
    0& 1& 2& 3& 4& 5& 6& 7& 8& 9&10&11\\
    1& 0& 4& 5& 2& 3& 7& 6&10&11& 8& 9\\
    2& 5& 0& 4& 3& 1& 8&11& 6&10& 9& 7\\
    3& 4& 5& 0& 1& 2& 9&10&11& 6& 7& 8\\
    4& 3& 1& 2& 5& 0&10& 9& 7& 8&11& 6\\
    5& 2& 3& 1& 0& 4&11& 8& 9& 7& 6&10\\
    6& 8& 7&11& 9&10& 0& 2& 1& 4& 5& 3\\
    7&10& 6& 9&11& 8& 1& 4& 0& 2& 3& 5\\
    8& 6& 9&10& 7&11& 2& 0& 5& 3& 1& 4\\
    9&11& 8& 7&10& 6& 3& 5& 4& 1& 2& 0\\
    10& 7&11& 8& 6 &9 &4 &1 &3&5& 0& 2\\
    11& 9&10& 6& 8& 7& 5& 3& 2& 0 &4&1
\end{array}
\end{displaymath}
Check that 1 is in the nucleus and $3^{-1}\cdot (1\cdot 3) = 2$. But $2$ is
not in the nucleus, since $4\cdot (6 \cdot 2) \ne (4 \cdot 6) \cdot 2$.
\end{example}

Fortunately, all is well for C-loops. We will use the following notation in the
proof of Proposition \ref{Pr:N}. Any element $x\in L$ determines two
permutations of $L$: the \emph{left translation} $L_x$ defined by $L_x(y)=xy$,
and the \emph{right translation} $R_x$ defined by $R_x(y)=yx$.

\begin{proposition}\label{Pr:N} The nucleus of a C-loop is a normal subgroup.
\end{proposition}
\begin{proof}
Let $N$ be the nucleus of a C-loop $L$. Our task is to show that $xN=Nx$ for
every $x\in L$, or, equivalently, that $x^{-1}nx\in N$ for every $x\in L$,
$n\in N$. Since the nuclei of a C-loop coincide, it suffices to show
$x^{-1}nx\in N_\lambda$, which in the language of translations becomes
$L_{x^{-1}nx}L_y = L_{(x^{-1}nx)y}$ for every $y\in L$.

Because squares of elements in a C-loop are in the nucleus and because
$x^2x^{-1}=x$, the last identity is equivalent to $L_{xnx}L_y = L_{(xnx)y}$,
which is what we prove below.

The following permutations coincide: $L_{xnx}$, $L_{n^{-1}(nx)^2}$ (by the left
inverse property and the right alternative property), $L_{n^{-1}}L_{(nx)^2}$
(since $(nx)^2\in N$), $L_{n^{-1}}L_nL_xL_nL_x$ (since $n\in N$), $L_xL_nL_x$
(by the left inverse property).

Using similar arguments, we see that $L_{x(n(xy))} = L_{(xn)^2(n^{-1}y)} =
L_{xn}L_{xn}L_{n^{-1}y} = L_xL_nL_xL_nL_{n^{-1}y} = L_xL_nL_xL_y$.

Therefore $L_{xnx}L_y = L_xL_nL_xL_y = L_{x(n(xy))}$. The last translation
$L_{x(n(xy))}$ is equal to $L_{(xnx)y}$, because $L_{xnx}=L_xL_nL_x$, and we
are done.
\end{proof}

\begin{proposition}\label{Pr:LN}
Let $L$ be a C-loop with nucleus $N$. Then $L/N$ is a Steiner loop.
\end{proposition}
\begin{proof}
We have $x^2\in N$ for every $x\in L$. Thus $L/N$ is an inverse property loop
of exponent $2$. By Lemma \ref{Lm:Folklore}, $L/N$ is a Steiner loop.
\end{proof}

The following Lemma will be useful in the next section.

\begin{lemma}\label{Lm:N2} There is no C-loop with nucleus of index $2$.
\end{lemma}
\begin{proof} Assume, for a contradiction, that $L$ is a C-loop with nucleus
$N$ of index $2$. Let $N$, $xN$ be the two cosets of $L/N$. We show that $x\in
N$.

Since the three nuclei of $L$ coincide, it suffices to show that $(ax)b=a(xb)$
for every $a$, $b\in L$. In fact, it suffices to prove this for all elements
$a$, $b\in xN=Nx$, since all other elements are nuclear. Let us write $a=cx$,
$b=xd$, for some $c$, $d\in N$. Since $c$, $d$, $x^2$ and $x^2d$ are all
nuclear, and since $x^2x = xx^2$, we have $(cx\cdot x)(xd) = (cx^2)(xd) =
c(x^2xd)=c(xx^2d) = (cx)(x^2d)=(cx)(x\cdot xd)$.
\end{proof}

\section{Admissible orders and the four smallest nonassociative C-loops}

\noindent We now construct the $4$ smallest nonassociative C-loops (1 of
order 10, 1 of order 12, and 2 of order 14). Three of these loops are
well-known Steiner loops. The remaining C-loop of order 12 belongs to an
infinite family of nonassociative noncommutative C-loops constructed here for
the first time. Although the four loops are constructed by hand, we do not
have sufficiently strong theoretical tools to show that no other
nonassociative C-loops of order less than $15$ exist. This is easily verified
by Mace4, though.

\subsection{Admissible orders}

Recall that Steiner loops of order $2$, $4$ and $8$ are elementary abelian
$2$-groups \cite{CR}.

\begin{proposition}\label{Pr:Orders} Let $L$ be a nonassociative C-loop of order
$n$ with nucleus $N$ of order $m$. Then
\begin{enumerate}
    \item[(i)] $n/m\equiv 2\pmod 6$ or $n/m\equiv 4\pmod 6$,
    \item[(ii)] $n$ is even,
    \item[(iii)] if $n=p^k$ for some prime $p$ and positive integer $k$, then
    $p=2$ and $k>3$.
\end{enumerate}
Moreover, there is a nonassociative non-Steiner C-loop of order $2^k$ for every
$k>3$.
\end{proposition}
\begin{proof}
Part (i) follows from Proposition \ref{Pr:LN} and from the already mentioned
fact that Steiner quasigroups of order $r$ exist if and only if $r\equiv 1\pmod
6$ or $r\equiv 3\pmod 6$. Part (ii) follows immediately from part (i).

Assume that $n=p^k$, $p$ a prime. By (ii), $p=2$. When $k<3$, $L$ must be a
group, since there is no nonassociative loop of order less than $5$.

Assume that $k=3$. If $m=1$, Proposition \ref{Pr:LN} implies that $L$ is a
Steiner loop of order $8$, thus the elementary abelian $2$-group of order
$8$. If $m=4$, we reach a contradiction by Lemma \ref{Lm:N2}. We were not
able to find a one-line argument that shows that there is no nonassociative
C-loop of order $8$ with nucleus of size $2$. It can be checked tediously by
hand.

Example \ref{Ex:16} gives a nonassociative non-Steiner
C-loop of order $16$. Direct products of this loop with
$2$-groups provide all needed examples.
\end{proof}

Let $L$, $n$, $m$ be as assumed in Proposition \ref{Pr:Orders}. The only
admissible values of $(n,m)$ with $n\le 14$ are then $(6,3)$, $(10,1)$,
$(10,5)$, $(12,3)$, $(12,6)$, $(14,1)$ and $(14,7)$. Lemma \ref{Lm:N2} further
reduces the possibilities to $(10,1)$, $(12,3)$ and $(14,1)$. As we shall see,
there is at least one nonassociative C-loop for each of these parameters.

\subsection{The smallest C-loop}
The smallest nonassociative commutative inverse property loop is of order $10$,
and it is unique. Its multiplication table is in Example \ref{Ex:10}. We can
see immediately that this loop has exponent $2$. It is therefore a Steiner
loop. By Lemma \ref{Lm:CS}, it is a nonassociative C-loop, hence the smallest
nonassociative C-loop.

\begin{example}\label{Ex:10}
The smallest nonassociative C-loop.
\begin{displaymath}
\begin{array}{cccccccccc}
    0& 1& 2& 3& 4& 5& 6& 7& 8& 9\\
    1& 0& 3& 2& 5& 4& 9& 8& 7& 6\\
    2& 3& 0& 1& 6& 8& 4& 9& 5& 7\\
    3& 2& 1& 0& 7& 9& 8& 4& 6& 5\\
    4& 5& 6& 7& 0& 1& 2& 3& 9& 8\\
    5& 4& 8& 9& 1& 0& 7& 6& 2& 3\\
    6& 9& 4& 8& 2& 7& 0& 5& 3& 1\\
    7& 8& 9& 4& 3& 6& 5& 0& 1& 2\\
    8& 7& 5& 6& 9& 2& 3& 1& 0& 4\\
    9& 6& 7& 5& 8& 3& 1& 2& 4& 0
\end{array}
\end{displaymath}
\end{example}

\subsection{The smallest noncommutative C-loop}
We now construct an infinite family of nonassociative noncommutative C-loops
whose smallest member is the smallest nonassociative noncommutative C-loop. The
construction is best approached via extensions of loops. Our notation is based
on \cite{DV}.

Let $G$ be a multiplicative group with neutral element $1$, and $A$ an abelian
group written additively with neutral element $0$. Any map $\mu:G\times G\to A$
satisfying $\mu(1,g)=\mu(g,1)=0$ for every $g\in G$ is called a \emph{factor
set}. When $\mu:G\times G\to A$ is a factor set, we can define multiplication
on $G\times A$ by
\begin{equation}\label{Eq:Mult}
    (g,a)(h,b)=(gh,a+b+\mu(g,h)).
\end{equation}
The resulting groupoid is clearly a loop with neutral element $(1,0)$. It will
be denoted by $(G,A,\mu)$. Additional properties of $(G,A,\mu)$ can be enforced
by additional requirements on $\mu$.

\begin{lemma} Let $\mu:G\times G\to A$ be a factor set. Then $(G,A,\mu)$ is a
C-loop if and only if
\begin{equation}\label{Eq:Ex}
    \mu(h,k)+\mu(h,hk)+\mu(g,h\cdot hk) = \mu(g,h)+\mu(gh,h) +\mu(gh\cdot h,k)
\end{equation}
for every $g$, $h$, $k\in G$.
\end{lemma}
\begin{proof}
The loop $(G,A,\mu)$ is a C-loop if and only if
\begin{displaymath}
    (g,a)((h,b)\cdot(h,b)(k,c)) = ((g,a)(h,b)\cdot(h,b))(k,c)
\end{displaymath}
holds for every $g$, $h$, $k\in G$ and every $a$, $b$, $c\in A$.
Straightforward calculation with $(\ref{Eq:Mult})$ shows that this happens if
and only if $(\ref{Eq:Ex})$ is satisfied.
\end{proof}

We call a factor set $\mu$ satisfying $(\ref{Eq:Ex})$ a \emph{C-factor set}.

When $G$ is an elementary abelian $2$-group, the equation $(\ref{Eq:Ex})$
reduces to
\begin{equation}\label{Eq:Ex2}
    \mu(h,k)+\mu(h,hk) = \mu(g,h)+\mu(gh,h).
\end{equation}
We now use a particular C-factor set to construct the above-mentioned family of
C-loops.

\begin{proposition}\label{Pr:Constr}
Let $n>2$ be an integer. Let $A$ be an abelian group of order $n$, and
$\alpha\in A$ an element of order bigger than $2$. Let $G=\{1,u,v,w\}$ be the
Klein group with neutral element $1$. Define $\mu:G\times G\to A$ by
\begin{displaymath}
    \mu(x,y) = \left\{
    \begin{array}{rl}
        \alpha,&\text{if }(x,y)=(v,w),\,(w,u),\,(w,w),\\
        -\alpha,&\text{if }(x,y)=(v,u),\\
        0,&\text{otherwise.}
    \end{array}\right.
\end{displaymath}
Then $(G,A,\mu)$ is a non-flexible $($hence nonassociative$)$ noncommutative
C-loop with nucleus $N=\{(1,a);\;a\in A\}$.
\end{proposition}
\begin{proof}
The map $\mu$ is clearly a factor set. It can be depicted as follows:
\begin{displaymath}
    \begin{array}{c|cccc}
        \mu&1&u&v&w\\
        \hline
        1&0&0&0&0\\
        u&0&0&0&0\\
        v&0&-\alpha&0&\alpha\\
        w&0&\alpha&0&\alpha
    \end{array}
\end{displaymath}
To show that $C=(G,A,\mu)$ is a C-loop, we verify $(\ref{Eq:Ex2})$.

Since $\mu$ is a factor set, there is nothing to prove when $h=1$. Assume that
$h=u$. Then $(\ref{Eq:Ex2})$ becomes $\mu(u,k)+\mu(u,uk)=\mu(g,u)+\mu(gu,u)$,
and both sides of this equation are equal to $0$, no matter what $k$, $g\in G$
are. Assume that $h=v$. Then $(\ref{Eq:Ex2})$ becomes
$\mu(v,k)+\mu(v,vk)=\mu(g,v)+\mu(gv,v)$, and both sides of this equation are
again equal to $0$. Assume that $h=w$. Then $(\ref{Eq:Ex2})$ becomes
$\mu(w,k)+\mu(w,wk)=\mu(g,w)+\mu(gw,w)$, and both sides of this equation are
equal to $\alpha$.

Since $\alpha\ne 0$, the C-loop $C$ is not commutative. By Corollary
\ref{Cr:Nuc}, the three nuclei of $C$ coincide and will be denoted by $N$.
Let $a\in A$. Since $\alpha\ne -\alpha$, we have
$(u,a)(v,a)\cdot(u,a)=(v,3a+\alpha)\ne (v,3a-\alpha) = (u,a)\cdot(v,a)(u,a)$.
This shows that: (i) $C$ is not flexible, (ii) $(u,a)$, $(v,a)\not\in N$.
Similarly, $(u,a)(w,a)\cdot (u,a)\ne (u,a)\cdot (w,a)(u,a)$ shows that
$(w,a)\not\in N$. Finally, for $g$, $h\in G$ and $b$, $c\in A$ we have
$(1,a)(g,b)\cdot (h,c) = (gh, a+b+c+\mu(g,h)) = (1,a)\cdot(g,b)(h,c)$, and
$(1,a)\in N$ follows.
\end{proof}

\begin{corollary}\label{Cr:NNNuc} For any integer $n>1$ there is a nonassociative
noncommutative $C$-loop with nucleus of size $n$.
\end{corollary}
\begin{proof}
It remains to show that there is a nonassociative noncommutative $C$-loop with
nucleus of size $2$. Consider the octonion loop $O$ of order $16$. This loop is
Moufang. Recall that extra loops are precisely Moufang loops with squares in
the nucleus. Since the squares in $O$ are equal to $1$ or $-1$, $O$ is an extra
loop, hence a $C$-loop. It is well-known that $N(O)=\{1,-1\}$.
\end{proof}

\begin{remark} The bound $n>1$ of Corollary $\ref{Cr:NNNuc}$ cannot be improved,
since a C-loop with nucleus of size $1$ is Steiner, hence commutative, by
Proposition $\ref{Pr:LN}$.
\end{remark}

\begin{example}\label{Ex:12}
The smallest group $A$ satisfying the assumptions of Proposition
\ref{Pr:Constr} is the $3$-element cyclic group $\{0,1,2\}$. The construction
of Proposition \ref{Pr:Constr} with $\alpha=2$ then gives rise to the
smallest noncommutative nonassociative C-loop:
\begin{displaymath}
\begin{array}{cccccccccccc}
         0& 1& 2& 3& 4& 5& 6& 7& 8& 9&10&11\\
         1& 2& 0& 4& 5& 3& 7& 8& 6&10&11& 9\\
         2& 0& 1& 5& 3& 4& 8& 6& 7&11& 9&10\\
         3& 4& 5& 0& 1& 2& 9&10&11& 6& 7& 8\\
         4& 5& 3& 1& 2& 0&10&11& 9& 7& 8& 6\\
         5& 3& 4& 2& 0& 1&11& 9&10& 8& 6& 7\\
         6& 7& 8&10&11& 9& 0& 1& 2& 5& 3& 4\\
         7& 8& 6&11& 9&10& 1& 2& 0& 3& 4& 5\\
         8& 6& 7& 9&10&11& 2& 0& 1& 4& 5& 3\\
         9&10&11& 8& 6& 7& 3& 4& 5& 2& 0& 1\\
        10&11& 9& 6& 7& 8& 4& 5& 3& 0& 1& 2\\
        11& 9&10& 7& 8& 6& 5& 3& 4& 1& 2& 0
\end{array}
\end{displaymath}
The following properties of this loop will be revoked later: The associator
$2=[11,8,9]$ has order $3$. Note that $3\cdot 3=0$, $6\cdot 6=0$, $3\cdot
6=9$, but $9\cdot 9=2\ne 0$.
\end{example}

\subsection{C-loops of order $14$}
There are two nonisomorphic nonassociative C-loops of order 14, both of them
Steiner loops. These loops are well known. We include their multiplication
tables for the sake of completeness, and also because we will refer to some of
their properties later.

Their multiplication tables are given in Examples \ref{Ex:14a} and
\ref{Ex:14b}.

\begin{example}\label{Ex:14a}
One of the two nonassociative C-loops of order $14$.
\begin{displaymath}
\begin{array}{cccccccccccccc}
    0&1&2&3&4&5&6&7&8&9&10&11&12&13\\
    1&0&3&2&5&4&12&13&9&8&11&10&6&7\\
    2&3&0&1&6&7&4&5&11&12&13&8&9&10\\
    3&2&1&0&7&8&9&4&5&6&12&13&10&11\\
    4&5&6&7&0&1&2&3&10&13&8&12&11&9\\
    5&4&7&8&1&0&10&2&3&11&6&9&13&12\\
    6&12&4&9&2&10&0&11&13&3&5&7&1&8\\
    7&13&5&4&3&2&11&0&12&10&9&6&8&1\\
    8&9&11&5&10&3&13&12&0&1&4&2&7&6\\
    9&8&12&6&13&11&3&10&1&0&7&5&2&4\\
    10&11&13&12&8&6&5&9&4&7&0&1&3&2\\
    11&10&8&13&12&9&7&6&2&5&1&0&4&3\\
    12&6&9&10&11&13&1&8&7&2&3&4&0&5\\
    13&7&10&11&9&12&8&1&6&4&2&3&5&0
\end{array}
\end{displaymath}
As we know from Proposition \ref{Pr:Orders}, this loop has a trivial nucleus.
Thus the associator $10=[13,12,1]$ is not in the nucleus.
\end{example}

\begin{example}\label{Ex:14b}
One of the two nonassociative C-loops of order 14.
\begin{displaymath}
\begin{array}{cccccccccccccc}
0&1&2&3&4&5&6&7&8&9&10&11&12&13\\
1&0&3&2&5&4&11&12&13&10&9&6&7&8\\
2&3&0&1&6&7&4&5&11&12&13&8&9&10\\
3&2&1&0&7&8&9&4&5&6&12&13&10&11\\
4&5&6&7&0&1&2&3&10&13&8&12&11&9\\
5&4&7&8&1&0&10&2&3&11&6&9&13&12\\
6&11&4&9&2&10&0&13&12&3&5&1&8&7\\
7&12&5&4&3&2&13&0&9&8&11&10&1&6\\
8&13&11&5&10&3&12&9&0&7&4&2&6&1\\
9&10&12&6&13&11&3&8&7&0&1&5&2&4\\
10&9&13&12&8&6&5&11&4&1&0&7&3&2\\
11&6&8&13&12&9&1&10&2&5&7&0&4&3\\
12&7&9&10&11&13&8&1&6&2&3&4&0&5\\
13&8&10&11&9&12&7&6&1&4&2&3&5&0
\end{array}
\end{displaymath}
\end{example}

\subsection{The smallest nonassociative non-Steiner commutative C-loop}
We must go beyond $n=14$ to find a nonassociative non-Steiner commutative
C-loop. There is one of order $16$, and its multiplication table is in Example
\ref{Ex:16}. We were unable to determine the number of nonassociative C-loops
of order $16$.

\begin{example}\label{Ex:16}
A nonassociative non-Steiner commutative C-loop of the smallest possible order.
\begin{displaymath}
\begin{array}{cccccccccccccccc}
    0& 1& 2& 3& 4& 5& 6& 7& 8& 9&10&11&12&13&14&15\\
    1& 5& 6& 8& 0& 4&10& 2&11& 3& 7& 9&13&15&12&14\\
    2& 6& 0&12& 7&10& 1& 4&14&13& 5&15& 3& 9& 8&11\\
    3& 8&12& 0& 9&11&14&13& 1& 4&15& 5& 2& 7& 6&10\\
    4& 0& 7& 9& 5& 1& 2&10& 3&11& 6& 8&14&12&15&13\\
    5& 4&10&11& 1& 0& 7& 6& 9& 8& 2& 3&15&14&13&12\\
    6&10& 1&14& 2& 7& 5& 0&12&15& 4&13& 9& 3&11& 8\\
    7& 2& 4&13&10& 6& 0& 5&15&12& 1&14& 8&11& 3& 9\\
    8&11&14& 1& 3& 9&12&15& 5& 0&13& 4& 7& 2&10& 6\\
    9& 3&13& 4&11& 8&15&12& 0& 5&14& 1& 6&10& 2& 7\\
    10& 7& 5&15& 6& 2& 4& 1&13&14& 0&12&11& 8& 9& 3\\
    11& 9&15& 5& 8& 3&13&14& 4& 1&12& 0&10& 6& 7& 2\\
    12&13& 3& 2&14&15& 9& 8& 7& 6&11&10& 0& 1& 4& 5\\
    13&15& 9& 7&12&14& 3&11& 2&10& 8& 6& 1& 5& 0& 4\\
    14&12& 8& 6&15&13&11& 3&10& 2& 9& 7& 4& 0& 5& 1\\
    15&14&11&10&13&12& 8& 9& 6& 7& 3& 2& 5& 4& 1& 0
\end{array}
\end{displaymath}
\end{example}

\section{Power associativity, diassociativity, and Lagrange-like properties}

\noindent Many properties that we take for granted in groups do not hold in
C-loops. This section is concerned with subloops generated by one or two
elements; with the relations between the order of a loop, the order of a
subloop, and the order of an element; and with alike properties.

\subsection{Power associativity.}

A loop $L$ is \emph{power associative} if for every $x\in L$ and every $n\ge 0$
the power $x^n$ is well-defined.

Clearly, the powers $x^0$, $x$, and $x^2$ are always well-defined. Note that,
up to this point, we have carefully avoided all higher powers in our
calculations. But we did not have to:

\begin{proposition}[Fenyves] LC-loops are power associative.
\end{proposition}
\begin{proof}
The power $x^n$ is clearly well-defined for $n=0$, $1$, $2$, and since, by the
left alternative law, $xx^2=x^2x$, it is also well-defined for $n=3$.

Assume that $n>3$ and that $x^k$ is well-defined for every $k<n$. Let $r$,
$s>0$ be such that $r+s=n$. We now show that $x^rx^s$ can be rewritten
canonically as $x^{r+s-1}x$. Since $xx^s=x(xx^{s-1})=x^2x^{s-1}$, we can
assume that $r>1$. Then $x^rx^s = x(xx^{r-2})\cdot x^s$, which is by the
LC-identity and by the induction hypothesis equal to $(xx)(x^{r-2}x^s) =
(xx)(x^{r+s-3}x) = x(xx^{r+s-3})\cdot x = x^{r+s-1}x$.
\end{proof}

\begin{corollary} C-loops are power associative.
\end{corollary}

\subsection{Diassociativity}

\emph A loop $L$ is \emph{diassociative} if any subloop of $L$ generated by two
elements is a subgroup.

C-loops are not necessarily diassociative. To see this, consider the C-loop $L$
of order 12 with multiplication table given in Example \ref{Ex:12}. Note that
$\langle 5\rangle = \{0,1,2,3,4,5\}$. It is the nobvious from the
multiplication table that $\langle 5,6\rangle=L$. Thus $L$ is generated by two
elements, yet it is not associative.

In \cite{KKPARIF}, \emph{ARIF loops} are defined to be flexible loops
satisfying $(zx)(yxy) = (z(xyx))y$.

\begin{lemma}\label{Lm:ARIF} Flexible C-loops are ARIF loops.
\end{lemma}
\begin{proof} Since C-loops are alternative and have all squares in the
nucleus, we have $z(xy)^2 = (zx\cdot x^{-1})(xy)^2 = (zx)(x^{-1}(xy)^2) =
(zx)((x^{-1}\cdot xy)(xy)) = (zx)(yxy)$. We use this identity twice to obtain
$(zx)(yxy) = z(xy)^2 = z((xyx)x^{-1})^2 = (z(xyx))(x^{-1}(xyx)x^{-1}) =
(z(xyx))y$.
\end{proof}

\begin{lemma}\label{Lm:FD}
Flexible C-loops are diassociative. In particular, commutative C-loops are
diassociative.
\end{lemma}
\begin{proof}
By Lemma \ref{Lm:ARIF}, flexible C-loops are ARIF loops. By \cite{KKPARIF},
ARIF loops are diassociative.

C-loops are alternative. In the presence of commutativity, the two alternative
laws are not only equivalent to each other, but also to the flexible law.
\end{proof}

\subsection{Power alternativity}

Power alternativity is best expressed in terms of translations. A power
associative loop $L$ is \emph{left power alternative} if $L_{x^n}=L_x^n$ for
every $n>0$ and $x\in L$. Similarly, $L$ is \emph{right power alternative} if
$R_{x^n}=R_x^n$ for every $n>0$ and $x\in L$. Power associative loops that are
both left and right power alternative are called \emph{power alternative}.

\begin{lemma}\label{Lm:LCPA} LC-loops are left power alternative.
\end{lemma}
\begin{proof}
We have $L_{x^2}=L_x^2$ by left alternativity. Let $n>2$, and assume that
$L_{x^m}=L_x^m$ for every $m<n$. We have $L_{x^n}(y)=x^ny=(x^{n-2}x^2)y$ by
power associativity. Since LC-loops are middle nuclear square, we have
$L_{x^n}(y)=x^{n-2}(x^2y)$, which, by induction, is equal to $L_x^n(y)$.
\end{proof}

\begin{corollary} C-loops are power alternative.
\end{corollary}

\subsection{Lagrange-like properties}

A finite loop $L$ is said to have the \emph{weak Lagrange property} if the
order of any subloop of $L$ divides the order of $L$. A finite loop $L$ has the
\emph{weak monogenic Lagrange property} if the order of any monogenic subloop
of $L$ divides the order of $L$.

To any weak property, there is its strong version: A loop has the \emph{strong
property P} if every subloop of $L$ has the weak property $P$. The weak
property does not always imply its strong version. For instance, there are
loops with the weak but not the strong Lagrange property.

Steiner loops (and hence C-loops) do not have the weak Lagrange property. This
is illustrated by the loop in Example \ref{Ex:14a} that is of order 14 and
possesses a subloop $\{0,1,2,3\}$ of order $4$.

However, C-loops have the strong monogenic Lagrange property. We can establish
this by imitating the proof of the Lagrange theorem for groups.

\begin{lemma}\label{Lm:SMLP}
Let $L$ be a finite loop that is left power alternative and has the left
inverse property. Then $L$ has the strong monogenic Lagrange property.
\end{lemma}
\begin{proof}
It suffices to prove that $L$ has the weak monogenic Lagrange property.

Let $x\in L$, $H=\langle x\rangle$. We claim that any two right cosets of $H$
are either disjoint or coincide.

Let $Hy$, $Hz$ be two such cosets. Assume that $u\in Hy\cap Hz$. Then
$u=x^ny=x^mz$ for some $n$, $m\ge 0$. By the left inverse property,
$z=x^{-m}(x^ny)$. By the left power alternative law, $z=x^{n-m}y$. Then $Hz =
H(x^{n-m}y) = \{x^r(x^{n-m}y);\;r\ge 0\} = \{x^{r+n-m}y;\;r\ge 0\}=Hy$, by the
left power alternative law again.
\end{proof}

\begin{corollary}\label{Cr:SMLP}
Let $x$ be an element of a finite LC-loop $L$. Then the order of $x$ divides
the order of $L$.
\end{corollary}
\begin{proof}
LC-loops are left power alternative, by Lemma \ref{Lm:LCPA}, and have the left
inverse property, by Proposition \ref{Pr:Sum}(ii). We are done by Lemma
\ref{Lm:SMLP}.
\end{proof}

\subsection{Cauchy-like properties}

A finite power associative loop is said to have the \emph{weak Cauchy property}
if for any prime $p$ dividing the order of the loop there is an element of
order $p$.

Since there are Steiner loops of order different from $2^k$, yet all Steiner
loops have exponent $2$, it is clear that Steiner loops (and thus C-loops) do
not have the weak Cauchy property. Example \ref{Ex:10} illustrates this nicely
(and minimally) for $p=5$.

\subsection{$2$-loops}

Since our main structural result for commutative C-loops
(Corollary \ref{Cr:Decomposition}) requires the notion of a
$2$-C-loop, let us talk about $2$-loops.

A group $G$ of order $n$ is said to be a \emph{$2$-group} if $n=2^r$ for some
$r$, or, equivalently, if $G$ is of exponent $2^s$ for some $s$.

The trouble with power associative loops is that the above two properties are
not necessarily equivalent. The smallest possible counterexample is a
nonassociative power associative loop of order $5$ and exponent $2$, cf.
\cite[Example I.4.5]{Pf}.

Throughout this paper, we therefore postulate: A finite power associative loop
$L$ is said to be a \emph{$2$-loop} if it is of exponent $2^s$, for some $s$.

Note that if $L$ is a C-loop of order $2^k$ then $L$ is a $2$-loop, by
Corollary \ref{Cr:SMLP}. The converse is not true, even for the smaller class
of Steiner loops, as is demonstrated by the nonassociative Steiner loop of
order $10$.

\section{Square roots of unity}

\noindent Example \ref{Ex:12} demonstrates that the subset $K=\{x\in
L;\;x^2=e\}$ of a C-loop $L$ is not necessarily a subloop of $L$. We are
going to see that in the commutative case, $K$ is not only a subloop, but a
normal subloop.

\begin{lemma}\label{Lm:xyxy} Let $L$ be a commutative, alternative, inverse property
loop. Then $(xy)^2=x^2y^2$ for every $x$, $y\in L$.
\end{lemma}
\begin{proof} Consider $u=x^{-1}(x^{-1}\cdot(xy)(xy))$. By alternativity,
$u=x^{-1}(x^{-1}(xy)\cdot (xy))$. By the inverse property, $u=x^{-1}(y\cdot
xy)$. By commutativity, $u=x^{-1}(xy\cdot y)$. By alternativity,
$u=x^{-1}(x\cdot yy)$. Finally, by the inverse property, $u=yy$. Thus
$(xy)(xy)=x(xu)=x(x(yy))=(xx)(yy)$.
\end{proof}

\begin{proposition}\label{Pr:LK} Let $L$ be a commutative $C$-loop, and let $K$
consist of all elements of exponent $2$ in $L$. Then $K$ is a normal subloop of
$L$ and $L/K$ is a group.
\end{proposition}
\begin{proof}
By Lemma \ref{Lm:xyxy}, the map $x\mapsto x^2$ is an endomorphism of $L$. Its
kernel $K$ is therefore a normal subloop of $L$.

It remains to show that $L/K$ is associative. This is true if and only if
$((xy)z)^{-1}(x(yz)\cdot u)\in K$ for every $x$, $y$, $z\in L$, $u\in K$, or,
equivalently, if $((xy)z)^2=(x(yz)\cdot u)^2$. Since squaring is a homomorphism
and since all squares are in the nucleus, we can rewrite the last equation as
$x^2y^2z^2=x^2y^2z^2u^2$, which certainly holds, as $u\in K$.
\end{proof}

\begin{corollary} Let $L$ be a commutative C-loop and let $A$ be the subloop
of $L$ generated by all associators $[x,y,z]$, where $x$, $y$, $z\in L$. Then
$A$ is of exponent at most $2$.
\end{corollary}
\begin{proof}
Let $K$ be as in Proposition \ref{Pr:LK}. Since $L/K$ is associative, all
associators must be in $K$. Thus $A\subseteq K$.
\end{proof}

Even in the noncommutative case we can say something about the associators.

\begin{lemma} Let $L$ be a loop with normal nucleus. Then all associators
of $L$ commute with all nuclear elements of $L$.
\end{lemma}
\begin{proof}
This is \cite[Lemma 4.2(vii)]{KKPCC}.
\end{proof}

\begin{corollary} Let $L$ be a C-loop. Then all associators of $L$ commute with all
nuclear elements. In particular, associators commute with all squares.
\end{corollary}

It would be nice if products of associators were again associators.
Unfortunately, this fails already for extra loops (hence for C-loops), by
\cite{KKnew}.

\section{An analogy between extra loops and C-loops}

\noindent The smallest variety of nonassociative loops of Bol-Moufang type is
that of extra loops (cf.\ \cite{PhVo}). We would like to describe an analogy
between extra loops and commutative C-loops.

A loop $L$ is \emph{conjugacy closed}, if for every $x$, $y\in L$,
$L_x^{-1}L_yL_x$ is a left translation, and $R_x^{-1}R_yR_x$ is a right
translation.

It is well known that extra loops are exactly conjugacy closed Moufang loops
(see, for instance, \cite{KKPCC}). Basarab \cite{Ba} showed that $L/N$ is an
abelian group for any conjugacy closed loop $L$ and its nucleus $N$. (Also see
\cite{KKPCC}, \cite{Dr}.) Since extra loops are also precisely Moufang loops
with squares in the nucleus, $L/N$ is an elementary abelian $2$-group whenever
$L$ is extra. Proposition \ref{Pr:LN} shows that $L/N$ is a Steiner loop when
$L$ is a C-loop. Since Steiner loops are commutative of exponent $2$, they
differ from elementary abelian $2$-groups ``only'' in their lack of
associativity.

The analogy can be extended little further for commutative C-loops. We have
seen that all associators in a commutative C-loop are of order at most $2$. The
same is true for extra loops.

However, all associators of an extra loop are in the nucleus. This is not the
case for commutative C-loops, as Example \ref{Ex:14a} illustrates.

\section{Decomposition for finite commutative C-loops}

\noindent We finish this paper with a decomposition theorem for finite
commutative C-loops.

\begin{lemma}\label{Lm:UV}
Let $L$ be a finite commutative C-loop. Let $U=\{x\in L;\; |x|$ is a power of
$2\}$, $V=\{x\in L;\; |x|$ is relatively prime to $2\}$. Then:
\begin{enumerate}
\item[(i)] $U\le L$, $V\le L$,

\item[(ii)] $V$ is contained in the nucleus of $L$, hence $V$ is a commutative
group,

\item[(iii)] $V\unlhd L$,

\item[(iv)] $U\unlhd L$,

\item[(v)] $UV=\{uv;\;u\in U,\,v\in V\}=L$,

\item[(vi)] $U\cap V=\{e\}$.
\end{enumerate}
\end{lemma}
\begin{proof}
First of all, by commutativity and diassociativity (Lemma \ref{Lm:FD}), we have
$(xy)^n=x^ny^n$ for every $x$, $y\in L$ and every integer $n$.

Let $x$, $y\in U$. Let $n$ be the least common multiple of $|x|$, $|y|$. Since
$x$, $y$ are powers of $2$, $n$ is a power of $2$ (the maximum of $|x|$ and
$|y|$). As $(xy)^n=x^ny^n=e$, we see that $|xy|$ divides $n$, and is therefore
a power of $2$.

Let $x$, $y\in V$. Let $n$ be the least common multiple of $|x|$, $|y|$. Since
both $x$, $y$ are relatively prime to $2$, so is $n$. As $(xy)^n=x^ny^n=e$, we
see that $|xy|$ divides $n$, and is therefore relatively prime to $2$. We have
proved (i).

Let $x\in V$. We want to show that $x$ belongs to the nucleus of $L$. Let
$n=|x|$. Then $n+1=2m$ is even, and $x=x^{n+1}=(x^m)^2$ is a square. Since
C-loops are nuclear square loops, $x$ is in the nucleus.

Any subloop contained in the center of $L$ is normal in $L$. By (ii), $V$ is
contained in the center of $L$, and so $V\unlhd L$.

We now show that $U$ is normal in $L$. Thanks to commutativity, all we have to
show is that $x(yU)=(xy)U$ for every $x$, $y\in L$. This is equivalent to
showing that $z=(xy)^{-1}(x(yu))\in U$ for every $u\in U$. Let $s=2^k$ be the
order of $u$. Then $z^s = x^{-s}y^{-s}x^sy^su^s=(xy)^{-s}(xy)^s=e$, by Lemma
\ref{Lm:xyxy}. Thus the order of $z$ divides $s=2^k$, and $z\in U$ follows.

Since (vi) follows immediately from the definition of $U$ and $V$, it remains
to prove that $UV=L$. Let $x\in L$, and let $|x|=2^ks$, where $k>0$ and $s>0$
is an odd integer. (There is nothing to prove when $k=0$ or $s=0$.) Since
$2^k$, $s$ are relatively prime, there are integers $m$, $n$ such that
$1=m2^k+ns$ \cite[Theorem 2-4]{Bu}. Then $x=uv$, where $u=x^{ns}$,
$v=x^{2^km}$. Since $u^{(2^k)}=x^{ns2^k}=1$, we see that $|u|$ divides $2^k$,
hence $|u|$ is a power of $2$, and $u\in U$ follows. Similarly, since
$v^s=x^{2^kns}=1$, we see that $|v|$ divides $s$, hence $|v|$ is odd, and $v\in
V$ follows.
\end{proof}

Universal algebraists define loops equivalently as sets with three binary
operations $\cdot$, $\ldiv$, $\rdiv$, and one nullary operation $e$ such that
$x\cdot(x\ldiv y)=y$, $x\rdiv y\cdot y = x$, $(x\cdot y)\rdiv y = x$,
$x\ldiv(x\cdot y)=y$, $x\rdiv x=x\ldiv x=e$. Thus $x\ldiv y$ is the solution
$z$ to the equation $x\cdot z=y$, and similarly for $x\rdiv y$.

We will use this notation in the proof of the following theorem, that could be
called \emph{the internal direct product for loops}. The theorem appears in a
more general form in Bruck's book \cite[Lemma IV.5.1]{Br}. Since he does not
give a proof, we provide it.

\begin{theorem}[Bruck]\label{Th:DPL}
Let $L$ be a loop with normal subloops $K$, $H$ such that $K\cap H=\{e\}$,
$KH=\{kh;\;k\in K$, $h\in H\}=L$. Then $L$ is the direct product of $K$, $H$.
\end{theorem}
\begin{proof}
We first show that any element $x\in L$ decomposes uniquely as a product $kh$,
$k\in K$, $h\in H$. At least one decomposition exists since $KH=L$. Let
$k_0h_0=k_1h_1$ be two such decompositions. Then $k_0=(k_1h_1)/h_0$. Now,
$(k_1h_1)/h_0$ belongs to $(Kh_1)/(Kh_0) = K(h_1/h_0)$ (since $K$ is a normal
subloop), and there is therefore $k\in K$ such that $k_0 = k(h_1/h_0)$. Then
$k\ldiv k_0=h_1/h_0$ belongs to $K\cap H=\{e\}$, and thus $h_1=h_0$, which in
turn implies $k_0=(k_1h_1)/h_0=(k_1h_0)/h_0=k_1$.

Define $f:K\times H\to L$ by $(k,h)\mapsto kh$. By the preceding paragraph, $f$
is one-to-one and onto. It remains to show that $f$ is a homomorphism, i.e.,
that $(k_0h_0)(k_1h_1)=k_0k_1\cdot h_0h_1$ for every $k_0$, $k_1\in K$, $h_0$,
$h_1\in H$. Since $K$ is normal, we have $x=(k_0h_0)(k_1h_1)\in Kh_0\cdot Kh_1
= K(h_0h_1)$, and there is $k\in K$ such that $x=k(h_0h_1)$. Since $H$ is
normal, we have $x\in k_0H\cdot k_1H = (k_0k_1)H$, and there is $h\in H$ such
that $x=(k_0k_1)h$. Since $x$ has a unique decomposition, we must have
$k=k_0k_1$, and so $x=(k_0k_1)(h_0h_1)$.
\end{proof}

\begin{theorem}\label{Th:Decomposition}
Let $L$ be a finite commutative C-loop. Then $L=U\times V$, where $U=\{x\in
L;\;|x|$ is a power of $2\}$, $V=\{x\in L;\;|x|$ is odd$\}$.
\end{theorem}
\begin{proof}
Combine Lemma \ref{Lm:UV} and Theorem \ref{Th:DPL}.
\end{proof}

\begin{corollary}\label{Cr:Decomposition}
Every finite commutative C-loop is a direct product of a finite commutative
group and a finite commutative $2$-C-loop, and \emph{vice versa}.
\end{corollary}

It is worth noting that we did not assume finiteness of the
loop in this section, only the fact that every element has
finite order. Power associative loops with all elements of
finite order are called \emph{torsion loops}. Thus, Lemma
\ref{Lm:UV}, Theorem \ref{Th:Decomposition} and Corollary
\ref{Cr:Decomposition} remain valid if all occurrences of
``finite'' in their statements are replaced by ``torsion''.

We conclude this paper with yet another analogy between C-loops and Moufang
loops. Theorem \ref{Th:Decomposition} shows that finite commutative C-loops are
of the form $U\times V$, where $U$ consists of elements whose order is a power
of $2$, and $V$ consist of elements whose order is relatively prime to $2$. By
\cite[Corollary of Theorem 3]{Gl}, finite commutative Moufang loops are of the
from $U\times V$, where $U$ consists of elements whose order is a power of $3$,
and $V$ consist of elements whose order is relatively prime to $3$.

\section{Acknowledgement}

\noindent We thank Michael K. Kinyon, Indiana University South Bend, for many
useful comments, especially for the proof of Theorem \ref{Th:DPL}; Scott
Feller, Department of Chemistry, Wabash College, for generously allowing us
to use his Beowulf cluster for our Mace4 computations; the anonymous referee
who suggested a deeper look at C-loops of small order---we came up with the
construction of Proposition \ref{Pr:Constr} in response to his/her inquiry.


\bibliographystyle{plain}

\end{document}